\documentclass[11pt]{amsart}
\usepackage{amsmath,amsfonts,amsthm,mathrsfs}
\usepackage{amssymb}
\usepackage{graphics}
\usepackage[all]{xy}
                                                                           
\usepackage[unicode,bookmarks]{hyperref}
\usepackage[usenames,dvipsnames]{xcolor}
\hypersetup{colorlinks=true,citecolor=NavyBlue,linkcolor=Brown,urlcolor=Orange}
                                                                           
\usepackage[alphabetic,initials]{amsrefs}

 \usepackage{enumitem}

  \usepackage{chngcntr}
                                                                
 \usepackage{tikz}
                                                                
 \usepackage{tikz-cd}
                                                                
 \usetikzlibrary{matrix,arrows}
 \newlength{\myarrowsize} 
                                                                
\pgfarrowsdeclare{cmto}{cmto}{
\pgfsetdash{}{0pt} 
	\pgfsetbeveljoin 
\pgfsetroundcap 
\setlength{\myarrowsize}{0.6pt}
	\addtolength{\myarrowsize}{.5\pgflinewidth}
  \pgfarrowsleftextend{-4\myarrowsize-.5\pgflinewidth} 
	\pgfarrowsrightextend{.8\pgflinewidth}
  }{
 \setlength{\myarrowsize}{0.6pt} 
\addtolength{\myarrowsize}{.5\pgflinewidth}  
    \pgfsetlinewidth{0.5\pgflinewidth}
  \pgfsetroundjoin
\pgfpathmoveto{\pgfpoint{1.5\pgflinewidth}{0}}
\pgfpatharc{-109}{-170}{4\myarrowsize}
\pgfpatharc{10}{189}{0.58\pgflinewidth and 0.2\pgflinewidth}
\pgfpatharc{-170}{-115}{4\myarrowsize+\pgflinewidth}
\pgfpathclose
\pgfusepathqfillstroke
\pgfpathmoveto{\pgfpoint{1.5\pgflinewidth}{0}}
\pgfpatharc{109}{170}{4\myarrowsize}
\pgfpatharc{-10}{-189}{0.58\pgflinewidth and 0.2\pgflinewidth}
\pgfpatharc{170}{115}{4\myarrowsize+\pgflinewidth}
\pgfpathclose
 \pgfusepathqfillstroke
\pgfsetlinewidth{2\pgflinewidth}
}
                                                            
                                                            \pgfarrowsdeclare{cmonto}{cmonto}{
                                                            	\pgfsetdash{}{0pt} 
                                                            	\pgfsetbeveljoin 
                                                            	\pgfsetroundcap 
                                                            	\setlength{\myarrowsize}{0.6pt}
                                                            	\addtolength{\myarrowsize}{.5\pgflinewidth}
                                                            	\pgfarrowsleftextend{-4\myarrowsize-.5\pgflinewidth} 
                                                            	\pgfarrowsrightextend{.8\pgflinewidth}
                                                            }{
                                                            \setlength{\myarrowsize}{0.6pt} 
                                                            \addtolength{\myarrowsize}{.5\pgflinewidth}  
                                                            \pgfsetlinewidth{0.5\pgflinewidth}
                                                            \pgfsetroundjoin
                                                            \pgfpathmoveto{\pgfpoint{1.5\pgflinewidth}{0}}
                                                            \pgfpatharc{-109}{-170}{4\myarrowsize}
                                                            \pgfpatharc{10}{189}{0.58\pgflinewidth and 0.2\pgflinewidth}
                                                            \pgfpatharc{-170}{-115}{4\myarrowsize+\pgflinewidth}
                                                            \pgfpathclose
                                                            \pgfusepathqfillstroke
                                                            \pgfpathmoveto{\pgfpoint{1.5\pgflinewidth}{0}}
                                                            \pgfpatharc{109}{170}{4\myarrowsize}
                                                            \pgfpatharc{-10}{-189}{0.58\pgflinewidth and 0.2\pgflinewidth}
                                                            \pgfpatharc{170}{115}{4\myarrowsize+\pgflinewidth}
                                                            \pgfpathclose
                                                            \pgfusepathqfillstroke
                                                            \pgfpathmoveto{\pgfpoint{1.5\pgflinewidth-0.3em}{0}}
                                                            \pgfpatharc{-109}{-170}{4\myarrowsize}
                                                            \pgfpatharc{10}{189}{0.58\pgflinewidth and 0.2\pgflinewidth}
                                                            \pgfpatharc{-170}{-115}{4\myarrowsize+\pgflinewidth}
                                                            \pgfpathclose
                                                            \pgfusepathqfillstroke
                                                            \pgfpathmoveto{\pgfpoint{1.5\pgflinewidth-0.3em}{0}}
                                                            \pgfpatharc{109}{170}{4\myarrowsize}
                                                            \pgfpatharc{-10}{-189}{0.58\pgflinewidth and 0.2\pgflinewidth}
                                                            \pgfpatharc{170}{115}{4\myarrowsize+\pgflinewidth}
                                                            \pgfpathclose
                                                            \pgfusepathqfillstroke
                                                            \pgfsetlinewidth{2\pgflinewidth}
                                                        }
                                                        
                                                        \pgfarrowsdeclare{cmhook}{cmhook}{
                                                        	\pgfsetdash{}{0pt} 
                                                        	\pgfsetbeveljoin 
                                                        	\pgfsetroundcap 
                                                        	\setlength{\myarrowsize}{0.6pt}
                                                        	\addtolength{\myarrowsize}{.5\pgflinewidth}
                                                        	\pgfarrowsleftextend{-4\myarrowsize-.5\pgflinewidth} 
                                                        	\pgfarrowsrightextend{.8\pgflinewidth}
                                                        }{
                                                        \setlength{\myarrowsize}{0.6pt} 
                                                        \addtolength{\myarrowsize}{.5\pgflinewidth}  
                                                        \pgfsetdash{}{0pt}
                                                        \pgfsetroundcap
                                                        \pgfpathmoveto{\pgfqpoint{0pt}{-4.667\pgflinewidth}}
                                                        \pgfpathcurveto
                                                        {\pgfqpoint{4\pgflinewidth}{-4.667\pgflinewidth}}
                                                        {\pgfqpoint{4\pgflinewidth}{0pt}}
                                                        {\pgfpointorigin}
                                                        \pgfusepathqstroke
                                                    }
                                                    
                                                \newenvironment{diagram*}[2]{%
                                                	\[%
                                                	\begin{tikzpicture}[>=cmto,baseline=(current bounding box.center),%
                                                	to/.style={->,font=\scriptsize,cap=round},%
                                                	into/.style={cmhook->,font=\scriptsize,cap=round},%
                                                	onto/.style={-cmonto,font=\scriptsize,cap=round},%
                                                	math/.style={matrix of math nodes, row sep=#2, column sep=#1,%
                                                		text height=1.5ex, text depth=0.25ex}]%
                                                }{%
                                                \end{tikzpicture}%
                                                \]%
                                                \ignorespacesafterend%
                                            }
                                            
                                            \newcommand{\norm}[1]{\lVert#1\rVert}

                                            \newcommand{\abs}[1]{\lvert #1 \rvert}

                                            \def\overbar#1#2#3{{%
                                            		\setbox0=\hbox{\displaystyle{#1}}%
                                            		\dimen0=\wd0
                                            		\advance\dimen0 by -#2 
                                            		\vbox {\nointerlineskip \moveright #3 \vbox{\hrule height 0.3pt width \dimen0}%
                                            			\nointerlineskip \vskip 1.5pt \box0}%
                                            	}}
                                            	
                                            	\makeatletter
                                            	\let\@@seccntformat\@seccntformat
                                            	\renewcommand*{\@seccntformat}[1]{%
                                            		\expandafter\ifx\csname @seccntformat@#1\endcsname\relax
                                            		\expandafter\@@seccntformat
                                            		\else
                                            		\expandafter
                                            		\csname @seccntformat@#1\expandafter\endcsname
                                            		\fi
                                            		{#1}%
                                            	}
                                            	\newcommand*{\@seccntformat@subsection}[1]{%
                                            		\textbf{\csname the#1\endcsname.}
                                            	}
                                            	\makeatother
                                            	
                                            	\makeatletter
                                            	\let\@paragraph\paragraph
                                            	\renewcommand*{\paragraph}[1]{%
                                            		\vspace{0.3\baselineskip}%
                                            		\@paragraph{\textit{#1}}%
                                            	}
                                            	\makeatother
                                            	
                                            	\counterwithin{equation}{section}
                                            	\counterwithin{figure}{section}
                                            	
                                            	\newtheorem{theorem}[equation]{Theorem}
                                            	\newtheorem*{theorem*}{Theorem}
                                            	\newtheorem{lemma}[equation]{Lemma}
                                            	\newtheorem*{lemma*}{Lemma}

                                            	\newtheorem*{proposition*}{Proposition}

                                            	\theoremstyle{definition}
                                            	\newtheorem{definition}[equation]{Definition}
                                            	\newtheorem*{definition*}{Definition}
                                            	\newtheorem{remark}[equation]{Remark}

                                            	\newtheorem{example}[equation]{Example}
                                            	\newtheorem*{example*}{Example}
                                            	\newtheorem*{problem*}{Problem}

                                            	\theoremstyle{plain}

                                            	\newcommand{\theoremref}[1]{\hyperref[#1]{Theorem~\ref*{#1}}}
                                            	\newcommand{\lemmaref}[1]{\hyperref[#1]{Lemma~\ref*{#1}}}
                                            	\newcommand{\definitionref}[1]{\hyperref[#1]{Definition~\ref*{#1}}}
                                            	\newcommand{\propositionref}[1]{\hyperref[#1]{Proposition~\ref*{#1}}}
                                            	\newcommand{\conjectureref}[1]{\hyperref[#1]{Conjecture~\ref*{#1}}}
                                            	\newcommand{\corollaryref}[1]{\hyperref[#1]{Corollary~\ref*{#1}}}
                                            	\newcommand{\exampleref}[1]{\hyperref[#1]{Example~\ref*{#1}}}

                                            	\makeatletter
                                            	\let\old@caption\caption
                                            	\renewcommand*{\caption}[1]{%
                                            		\setcounter{figure}{\value{equation}}%
                                            		\stepcounter{equation}%
                                            		\old@caption{#1}\relax%
                                            	}
                                            	\makeatother
                                            	
                                            	\newcounter{intro}
                                            	
                                            	\newtheorem{intro-conjecture}[intro]{Conjecture}
                                            	\newtheorem{intro-corollary}[intro]{Corollary}
                                            	\newtheorem{intro-theorem}[intro]{Theorem}
                                            	
                                            	\newcommand{\parref}[1]{\hyperref[#1]{\S\ref*{#1}}}
                                            	
                                            	\makeatletter
                                            	\newcommand*\if@single[3]{%
                                            		\setbox0\hbox{${\mathaccent"0362{#1}}^H$}%
                                            		\setbox2\hbox{${\mathaccent"0362{\kern0pt#1}}^H$}%
                                            		\ifdim\ht0=\ht2 #3\else #2\fi
                                            	}
                                            	\newcommand*\rel@kern[1]{\kern#1\dimexpr\macc@kerna}
                                            	\newcommand*\widebar[1]{\@ifnextchar^{{\wide@bar{#1}{0}}}{\wide@bar{#1}{1}}}
                                            	\newcommand*\wide@bar[2]{\if@single{#1}{\wide@bar@{#1}{#2}{1}}{\wide@bar@{#1}{#2}{2}}}
                                            	\newcommand*\wide@bar@[3]{%
                                            		\begingroup
                                            		\def\mathaccent##1##2{%
                                            			\if#32 \let\macc@nucleus\first@char \fi
                                            			\setbox\z@\hbox{$\macc@style{\macc@nucleus}_{}$}%
                                            			\setbox\tw@\hbox{$\macc@style{\macc@nucleus}{}_{}$}%
                                            			\dimen@\wd\tw@
                                            			\advance\dimen@-\wd\z@
                                            			\divide\dimen@ 3
                                            			\@tempdima\wd\tw@
                                            			\advance\@tempdima-\scriptspace
                                            			\divide\@tempdima 10
                                            			\advance\dimen@-\@tempdima
                                            			\ifdim\dimen@>\z@ \dimen@0pt\fi
                                            			\rel@kern{0.6}\kern-\dimen@
                                            			\if#31
                                            			\overline{\rel@kern{-0.6}\kern\dimen@\macc@nucleus\rel@kern{0.4}\kern\dimen@}%
                                            			\advance\dimen@0.4\dimexpr\macc@kerna
                                            			\let\final@kern#2%
                                            			\ifdim\dimen@<\z@ \let\final@kern1\fi
                                            			\if\final@kern1 \kern-\dimen@\fi
                                            			\else
                                            			\overline{\rel@kern{-0.6}\kern\dimen@#1}%
                                            			\fi
                                            		}%
                                            		\macc@depth\@ne
                                            		\let\math@bgroup\@empty \let\math@egroup\macc@set@skewchar
                                            		\mathsurround\z@ \frozen@everymath{\mathgroup\macc@group\relax}%
                                            		\macc@set@skewchar\relax
                                            		\let\mathaccentV\macc@nested@a
                                            		\if#31
                                            		\macc@nested@a\relax111{#1}%
                                            		\else
                                            		\def\gobble@till@marker##1\endmarker{}%
                                            		\futurelet\first@char\gobble@till@marker#1\endmarker
                                            		\ifcat\noexpand\first@char A\else
                                            		\def\first@char{}%
                                            		\fi
                                            		\macc@nested@a\relax111{\first@char}%
                                            		\fi
                                            		\endgroup
                                            	}
                                            	\makeatother

                                            	\usepackage{amsmath,amsfonts,mathrsfs}

                                            	\newcommand{\Pd}{\mathcal{P}}
                                            	\newcommand{\Om}{\mathcal{O}}

                                            	\def\F{{\mathcal F}}

                                            \newcommand{\J}{\mathcal{J}}

                                            \newcommand{\Di}{\textup{dim }}

\begin{document}
                                            	
                                            	\title{On 3-canonical maps of varieties of Albanese fiber dimension one}
                                            	
                                            		\author[Y. Chen]{Yuesen Chen}
                                            	\address{Department of Mathematics, Fudan University ,
                                            		Shanghai, China}
                                            	\email{{14110180003@fudan.edu.cn}}

                                            	\begin{abstract}
                                            		In the present paper, we study the (twisted) 3-canonical map of varieties of Albanese fiber dimension one. Based on a theorem about the regularity of direct image of canonical sheaves, we prove that the 3-canonical map is generically birational when the genus of a general fiber of the Albanese map is 2. 
                                            	\end{abstract}
                                            	
                                       	\maketitle
                                            	
                                            	\section{Introduction}
                                            	
                                            	Pluricanonical maps are an essential tool for understanding varieties of general type in birational geometry. As for irregular varieties, the most widely studied class is the varieties of maximal Albanese dimension (m.A.d. for short), whose pluricanonical maps have been researched extensively: Chen and Hacon \cite{ch} showed that the 3-canonical map of general type m.A.d. with $\chi(X,\omega_X)>0$ is birational, Jiang, Lahoz and Tirabassi \cite{jlt} then improved this result by eliminating the last condition. Barja et al. \cite{blnp} studied the bicanonical map of m.A.d. of general type under some natural conditions. Lahoz \cite{lah} and Zhang \cite{lzhang} also proved some criteria for the birationality of bicanonical map from different points of view.
                                            	
                                            	However, for a variety $X$ of general type with Albanese fiber dimension 1, the properties of the pluricanonical map has not yet been well-established. The most updated literature is \cite{js}, in which Jiang and Sun proved that $\abs{4K_X}$ is birational. It is then natural and interesting to study $\abs{3K_X}$ for varieties of Albanese fiber dimension 1 carefully. One might conjecture that $\abs{3K_X}$ is birational since $\abs{3K_C}$ is birational for a curve $C$ of genus $g\ge2$. Even though we are unable to prove (or disprove) the conjecture, we provide some evidence towards a positive answer to this conjecture:
                                            	
                                            	\begin{theorem}\label{main}
                                            		Let $X$ be a smooth complex projective variety of general type and of Albanese fiber dimension one, if the genus of a general fiber of the Albanese map is 2, then for general $P\in \mathrm{Pic}^0(X)$, $\abs{3K_X+P}$ is birational.
                                            	\end{theorem}
                                            	
                                            	We give a summary of the proof. By a standard argument, it suffices to prove the theorem for a general $P\in \mathrm{Pic}^0(X)_{tors}$. After constructing an intrinsic morphism, our proof is then divided into two parts: first we prove that the 3-canonical linear system separates two general fibers of the morphism; then we show that the 3-canonical linear system induces a birational map on each general fiber. Each step is heavily dependent on a decomposition theorem by Pareschi, Popa and Schnell (Theorem \ref{dt}).
                                            	
                                            	\bigskip
                                            	\section{Preliminaries}
                                            	
                                            	We work over $\mathbb{C}$, the field of complex numbers, throughout this paper.
                                            	
                                            	Let $X$ be a smooth projective variety and $a_X:X\rightarrow A_X$ the Albanese map. We say that $X$ is of Albanese fiber dimension $l$ if $\Di X-\Di a_X(X)=l$. In particular, $X$ is of Albanese fiber dimension 0 if and only if $X$ is of maximal Albanese dimension.
                                            	
                                            	Given an abelian variety $A$, let $\hat{A}=\mathrm{Pic}^0(A)$ be its dual, whose elements are line bundles on $A$ with trivial algebraic equivalence class. In some situation by an element of $\mathrm{Pic}^0(A)$ we also mean a divisor representing the line bundle.
                                            	
                                            	For a morphism $t:X\rightarrow A$ to an abelian variety and a coherent sheaf $\F$ on $X$, we denote by $V^i(\F,t)$ the $i$-th cohomological support loci:
                                            	\begin{equation*}
                                            	\{P\in \mathrm{Pic}^0(A)\mid H^i(X,\F\otimes \mathop{t^{*}}P)\ne 0\}.
                                            	\end{equation*}
                                            	In particular, if $t=a_X$, we will simply denote $V^i(\F,a_X)$ by $V^i(\F)$.
                                            	
                                            	We recall the following definitions due to Pareschi, Popa \cite{pp1} and Mukai. 
                                            	
                                            	\begin{definition}\label{GV}
                                            		Let $\F$ be a coherent sheaf on a smooth projective variety $X$.
                                            		\begin{enumerate}
                                            			\item The sheaf $\F$ is said to be a GV-sheaf if $\mathrm{codim}V^i(\F) \ge i$ for all $i \ge 0$.
                                            			\item If $X$ is an abelian variety, a non-zero sheaf $\F$ is said to be M-regular if $\mathrm{codim} V^i(\F) > i$ for every $i>0$. Furthermore, a M-regular sheaf $\F$ is said to be IT(0) if $V^i(\F)=\emptyset$ for all $i>0$. (We set $\dim \emptyset=-\infty$.)
                                            		\end{enumerate}
                                            	\end{definition}

                                            	In the next, we introduce the notion of continuously globally generation:
                                            	
                                            	\begin{definition}\label{cgg}Let $a:X\rightarrow A$ be a morphism from a projective variety to an abelian variety, $\F$ a coherent sheaf on $X$ and $T$ a subset of $ \hat{A}$. Then we have the continuous evaluation map
                                            		\begin{equation*}
                                            		ev_{T,\F}:\bigoplus\limits_{\alpha\in T}H^0(X,\F\otimes \mathop{a^*}\alpha^{-1})\otimes \mathop{a^*}\alpha\rightarrow \F.
                                            		\end{equation*}
                                            		Let $p$ be a closed point of $X$. The sheaf $\F$ is said to be \textit{continuously globally generated} at $p$ (CGG at $p$ for short), with respect to the morphism $a$, if the map $ev_{U,\F}$ is surjective at $p$, for all Zariski open dense subsets $U\subseteq\hat{A}$. Moreover, we will say that a sheaf is
                                            		CGG, when it is CGG for all $p$.
                                            	\end{definition}
                                            	
                                            	M-regular sheaves have the CGG property:
                                            	\begin{theorem}[\cite{blnp}]\label{Mtocgg}
                                            		Any M-regular sheaf on an abelian variety $A$ is continuously globally generated. If $\F$ is a M-regular sheaf on $A$, then $h^0(\F\otimes P)> 0$ for any $P\in\hat{A}$; if $\F$ is IT(0) then $h^0(\F\otimes P)$ is independent of the choice of $P$.
                                            	\end{theorem}

                                            	Pareschi, Popa and Schnell have proven a theorem that provides
                                            	certain M-regular properties for the pushforward of the canonical sheaf, which plays a critical role in this paper. This decomposition was originally discovered by Chen and Jiang \cite{cj} in the case of maximal Albanese dimension and then extended to irregular variety and higher direct images in general. We state the simplest form of their main theorem as below: 
                                            	
                                            	\begin{theorem}[\cite{pps}, Theorem A and \cite{lps}, Theorem C]\label{dt}
                                            		Let $f:X\rightarrow A$ be a morphism from a smooth projective variety $X$ to an abelian variety $A$. Then there exist finitely many quotients of abelian varieties $p_{B_i}: A\rightarrow B_i$ with connected kernel, finitely many M-regular sheaves $\F^j_{B_i}$ on each $B_i$ and torsion line bundles $Q^j_{B_i} \in \mathrm{Pic}^0(A)
                                            		$ such that
                                            		\begin{equation}
                                            		\mathop{f_*}\omega_X=\bigoplus\limits_{\substack{p_{B_i}:A\rightarrow B_i,\\ i\in I }}\bigoplus\limits_{j} \mathop{p_{B_i}^*} \F^j_{B_i}\otimes Q^j_{B_i}.
                                            		\end{equation}
                                            	\end{theorem}
                                            	
                                            	The following theorem by Jiang and Sun makes our conditions accessible:\medskip
                                            	\begin{theorem}[\cite{js}, Theorem 3.3 and Lemma 4.1]\label{js1}
                                            		If $X$ is a smooth projective variety of general type and of Albanese fiber dimension one, with Albanese morphism $a_X:X\rightarrow A$. Then the translates through 0 of all irreducible components of $ V^0(\omega_X,a_X)$ generate $\mathrm{Pic}^0(A)$.
                                            	\end{theorem}
                                            	
                                            	Finally, we will use the following results regarding M-regularity.
                                            	
                                            	\begin{theorem}\label{popas}Let $f:X\rightarrow A$ be a morphism from a smooth projective variety to an abelian variety, $Q\in\mathrm{Pic}^0(X)$ is a torsion element, $k$ an integer, then
                                            		\begin{enumerate}
                                            			\item If $k>0$, $\mathop{f_*}(\omega_X^k\otimes Q)$is GV. 
                                            			\item If $X$ is of general type, $k>1$, and $\mathop{f_*}(\omega_X^k\otimes Q)$ is non-zero, then $\mathop{f_*}(\omega_X^k\otimes Q)$ is IT(0).
                                            		\end{enumerate}
                                            		
                                            	\end{theorem}
                                            	\begin{proof}
                                            		For the proof of 1., see \cite{ps2014} and \cite[Theorem 2.2]{hp}. For the proof of 2., see \cite[Theorem 10.1 and Theorem 11.2]{hps},\cite[Proposition 2.3 and Corollary 11.2]{lps} and \cite[Lemma 2.2]{j1}.
                                            		
                                            	\end{proof}

                                            	\bigskip
                                            	
                                            	\section{Main Results}
                                            	Let $X$ be a smooth projective variety of general type and of Albanese fiber dimension one, with Albanese map $a:X\rightarrow A$.
                                            	Take its Stein factorization 	\begin{equation*}
                                            	\begin{xy}
                                            	(0,20)*+{X}="v1";(30,20)*+{ X_A}="v2";(15,0)*+{A}="v3";
                                            	{\ar@{->}_f"v1";"v2"};
                                            	{\ar@{->}_l"v2";"v3"};
                                            	{\ar@{->}_a"v1";"v3"}; 
                                            	\end{xy}
                                            	\end{equation*}
                                            	
                                            	We say $X$ is of \textit{Albanese fiber genus $m$} if the genus of a general fiber of $f$ is $\frac{m}{\textnormal{deg} l}$. 
                                            	
                                            	By the work of Jiang and Sun \cite[Remark 5.4]{js}, in this paper we only consider the case that the cohomological locus is not the whole Picard zero group i.e. $V^0(K_X)\neq \mathrm{Pic}^0(A)$, which is equivalent to the condition $\chi(\mathop{a_*}\omega_X)=0$. In fact, this situation does happen when $\dim X\ge 3$ (see Example \ref{examp}).

                                            	\begin{lemma}\label{tours}Notations as above, let $B_i$ be the quotient abelian varieties $p_{B_i}:A\rightarrow B_i$ in the application of Theorem \ref{dt}. Then those $\mathop{p_{B_i}^*}\mathrm{Pic}^0(B_i)$'s generate $\mathrm{Pic}^0(A)$ together.
                                            	\end{lemma}
                                            	\begin{proof}
                                            		By \cite[Lemma 3.3]{lps} we have $V^0(K_X)=\bigcup (\mathop{p_{B_i}^*}\textnormal{Pic}^0(B_i)-Q_{B_{i}}^j)$, hence the translates through 0 of all irreducible components of $ V^0(K_X)$ is just those $\mathop{p_{B_i}^*}\mathrm{Pic}^0(B_i)$'s, which generate $\mathrm{Pic}^0(A)$ by Theorem \ref{js1}. 
                                            	\end{proof}
                                            	\medskip
                                            	
                                            	Let $T_i$ be the subtori of $A$ corresponding to the kernel of the quotient maps $p_{B_i}:A\rightarrow B_i$ in Theorem \ref{dt},
                                            	then their intersection
                                            	\begin{equation*}
                                            	T:=\bigcap\limits_{i} T_i
                                            	\end{equation*} is a finite subgroup of $A$, by the lemma above.
                                            	
                                            	Denote by $d:X\rightarrow A/T$ the quotient map $A\rightarrow A/T$ composed with $a:X\rightarrow A$, and $f_{B_i}=p_{B_i}\circ a:X\rightarrow B_i$ the quotient map $A\rightarrow B_i$ composed with $a$. Because $A\rightarrow A/T$ is an isogeny, it's enough and much natural to deal with our problem over $d$ rather than $a$.
                                            	
                                            	\begin{equation*}
                                            	\begin{xy}
                                            	(20,20)*+{X}="v2";(20,0)*+{A}="v4";(35,0)*+{B_i}="v5";(20,-20)*+{A/T}="v3";
                                            	{\ar@{->}_{a}"v2";"v4"};
                                            	{\ar@{->}"v4";"v3"};
                                            	{\ar@{->}"v3";"v5"};
                                            	{\ar@/_1pc/_{d}"v2";"v3"};
                                            	{\ar@{->}_{p_{B_i}}"v4";"v5"};
                                            	{\ar@{->}^{f_{B_i}}"v2";"v5"};
                                            	\end{xy}
                                            	\end{equation*}
                                            	
                                            	Note that by Theorem \ref{js1} $V^0(K_X)\ne \emptyset$, hence $\mathop{a_*}\omega_X$ is non-zero; by the projection formula each term $\F_i, \ i\in I$ in Theorem \ref{dt} is torsion-free and of rank $\ge1$ on $f_{B_i}(X)$.
                                            	
                                            	We will use the notations above without any mention in this section.
                                            	
                                            	We have the following M-regularity:
                                            	
                                            	\begin{lemma}\label{genericM}
                                            		Notations as above, for any integer $k>1$ and any $P \in\mathrm{Pic}^0(X)_{tors}$ (the subgroup the torsion elements), $\mathop{f_{B_i*}}(\omega_X^k\otimes P)$ is IT(0) (hence M-regular) on $B_i$.
                                            	\end{lemma}
                                            	\begin{proof}
                                            		Since $V^0(2K_X)\ne \emptyset$ (Theorem \ref{js1}), $\mathop{a_*}\omega_X^2$ is non-zero therefore IT(0) (Theorem \ref{popas}, (2.)). Take any $P\in\mathrm{Pic}^0(X)_{tors}$, because $H^0(B_i, \mathop{f_{B_i*}}(\omega_X^2\otimes P))=H^0(A, \mathop{a_{*}}(\omega_X^2\otimes P))\ne 0$, $\mathop{f_{B_i*}}(\omega_X^k\otimes P)$ is a non-zero sheaf. Now apply Theorem \ref{popas}.
                                            	\end{proof}
                                            	
                                            	We are now devoted to prove our main result. Firstly we describe the 3-canonical linear systems over different fibers of $d$:
                                            	
                                            	\medskip
                                            	\begin{lemma}\label{df} For general $P\in \mathrm{Pic}^0(X)_{tors} $, and any $B_i$ appearing in Theorem \ref{dt}, the linear system $\arrowvert 3K_X+ P\arrowvert$ separates two general different points of $X$ lying in different fibers of $f_{B_i}$ over $f_{B_i}(X)$. Consequently, 	$\arrowvert 3K_X+ P\arrowvert$ separates two general different points lying in different general fibers of $d$ over $d(X)\subseteq A/T$.
                                            	\end{lemma}
                                            	\begin{proof}
                                            		Since $P\in \mathrm{Pic}^0(X)_{tors} $, $\mathop{f_{{B_i}*}}\mathcal{O}_X(2K_X+P)$ is M-regular by Lemma \ref{genericM}.
                                            		
                                            		Let $U_0\subseteq X$ be an open dense subset such that $f_{B_i}\vert_{U_0}$ and $a\vert_{U_0}$ are smooth, and that for any $b$ in the open dense subset $f_{B_i}(U_0) \subseteq f_{B_i}(X)$, we have the base change isomorphism $\mathop{f_{B_i*}}(\mathcal{O}_X(2K_X+P))\otimes \mathbb{C}_b\cong H^0(F_b,2K_{F_b}+P)$, where $F_b$ is the fiber of $f_{B_i}$ over $b$. Choose $D\in\abs{2K_X+P}$ an effective divisor on $X$ and denote $U_1:=U_0\cap D^c$ (here $D^c$ is the complement of $D$ in $X$).
                                            		
                                            		Let $x\in U_1$ be a closed point and $b_1=f_{B_i}(x)$. Theorem \ref{dt} implies that $p_{B_i}^{-1}(b_1)$ is contained in $a(X)\subseteq A$, therefore every irreducible component of the fiber $F_{b_1}$ over $b_1$ is a smooth variety of Albanese fiber dimension $\le 1$, and thus $H^0(F_{b_1},2K_{F_{b_1}}+P)$ is non-zero by M-regularity. Since $x$ is not a base point of $\abs{2K_{F_{b_1}}+P}$, we have an exact sequence (cf. \cite[Proposition 5.2]{js}) (here $\mathcal{I}_x$ is the ideal sheaf of $x$)
                                            		\begin{equation}\label{exa1}
                                            		0\rightarrow \mathop{f_{{B_i}*}}(\mathcal{O}_X(2K_X+P)\otimes \mathcal{I}_x)\rightarrow \mathop{f_{{B_i}*}}\mathcal{O}_X(2K_X+P)\rightarrow \mathbb{C}_{b_1}\rightarrow 0.
                                            		\end{equation}
                                            		
                                            		By employing a trick of Tirabassi \cite{t} (see \cite[Proposition 5.2]{js} for our case \footnote{In Jiang and Sun's original paper \cite{js}, they only considered the subsheaf  $\mathop{a_{*}}(\mathcal{O}_X(2K_X)\otimes\J(\norm{K_X}))$ (here $\J(\norm{K_X})$ denotes the asymptotic multiplier ideal of $K_X$). However, it's not hard to show that $\mathop{a_{*}}(\mathcal{O}_X(2K_X)\otimes\J(\norm{K_X}))=\mathop{a_{*}}(\mathcal{O}_X(2K_X))$ (cf. \cite[4.2]{jiangp}) and hence $\mathop{a_{*}}(\mathcal{I}_x\otimes\mathcal{O}_X(2K_X)\otimes\J(\norm{K_X}))=\mathop{a_{*}}(\mathcal{I}_x\otimes\mathcal{O}_X(2K_X))$, so their result applies.}), we know that for any $x$ in an open dense subset $U^\prime \subseteq X$, $\mathop{a_{*}}(\mathcal{I}_x\otimes\mathcal{O}_X(2K_X))$ is M-regular on $A$, and \begin{equation*}\{R\in\mathrm{Pic}^0(X)\vert \ x\in Bs(\abs{2K_X+R})\} \end{equation*}
                                            		has codimension $\ge 2$ in $\mathrm{Pic}^0(X)$. 
                                            		
                                            		Therefore, fix a general $P\in\mathrm{Pic}^0(X)_{tors}$, the closed subset
                                            		\begin{equation*}
                                            		\mathcal{B}(P):=\{R\in\hat{B_i}\vert \ x\in Bs(\abs{2K_X+P+R})\}
                                            		\end{equation*} has codimension $\ge 2$ in $\hat{B_i}$, for any choice of $x$ in a dense open subset $ U_2\subseteq U_1$ within $X$:
                                            		To see this, note that $V^1(\mathop{f_{B_i*}}(\mathcal{O}_X(2K_X+P)))=\emptyset$  by IT(0) property, and set
                                            		\begin{equation*}
                                            		\mathcal{D}=\{(x,Q)\in U^\prime\times \hat{B_i}\vert \ x \in Bs\abs{2K_X+P+Q}\},
                                            		\end{equation*} which is a closed subset in $ U^\prime\times \hat{B_i}$ of codimension $\ge2$ as above, by our generic choice of $P\in\mathrm{Pic}^0(X)_{tors}$. Let $p_1:X\times \hat{B_i}\rightarrow X$ be the first projection, then the required $U_2$ is nothing but the open subset of $U^\prime\cap U_1$ over which the fiber of $p_1\vert_{\mathcal{D}}$ has codimension $\ge2$ in $\hat{B_i}$.
                                            		
                                            		\medskip
                                            		
                                            		Take $x\in U_2$, since \begin{equation*}V^1(\mathop{f_{{B_i}*}}(\mathcal{O}_X(2K_X+P)\otimes \mathcal{I}_x))=\mathcal{B}(P)\cup V^1(\mathop{f_{{B_i}*}}\mathcal{O}_X(2K_X+P))\end{equation*} and \begin{equation*}V^i(\mathop{f_{{B_i}*}}(\mathcal{O}_X(2K_X+P)\otimes \mathcal{I}_x))= V^i(\mathop{f_{{B_i}*}}\mathcal{O}_X(2K_X+P)), \ i>1 ,\end{equation*} (by (\ref{exa1})) we conclude that $\mathop{f_{{B_i}*}}(\mathcal{O}_X(2K_X+P)\otimes \mathcal{I}_x)$ is M-regular by Definition \ref{GV}.
                                            		
                                            		\medskip
                                            		
                                            		For any $y\in U_2\subseteq X$ with $b_2:=f_{{B_i}}(y)\ne b_1=f_{{B_i}}(x)$, the map
                                            		\begin{equation}
                                            		\mathop{f_{{B_i}*}}(\mathcal{O}_X(2K_X+P)\otimes \mathcal{I}_x)\twoheadrightarrow \mathbb{C}_{b_2}
                                            		\end{equation} which is the pushforward of the natural surjection
                                            		\begin{equation*}
                                            		m_y: \mathcal{O}_X(2K_X+P)\otimes \mathcal{I}_x\rightarrow \mathbb{C}_y\cong \mathcal{O}_X(2K_X+P)\otimes \mathcal{I}_x\otimes \Om_y 	
                                            		\end{equation*} via $f_{{B_i}*}$, 	is again surjective, since $y$ is not a base point of $H^0(F_{b_2}, 2K_{F_{b2}}+P)\cong \mathop{f_{{B_i}*}}(\mathcal{O}_X(2K_X+P)\otimes \mathcal{I}_x)\otimes\mathbb{C}_{b_2}$. Then by the CGG property (Theorem \ref{Mtocgg}) of $\mathop{f_{{B_i}*}}(\mathcal{O}_X(2K_X+P)\otimes \mathcal{I}_x)$, there is an open dense subset $V\subseteq \hat{B_i}$, such that for any $R\in V$, the map
                                            		\begin{equation*}\label{firsts}
                                            		H^0(X,\mathcal{O}_X(2K_X+P+R)\otimes\mathcal{I}_x)  \xrightarrow{H^0(m_y)} \mathbb{C}_{y}
                                            		\end{equation*}is surjective.

                                            		On the other hand, denote by $Q$ the torsion line bundle as in Theorem \ref{dt} respect to a direct summand on ${B_i}$. Note that $\mathop{f_{B_i*}}\mathcal{O}_X(K_X-Q)$ is GV and $V^0(\mathop{f_{B_i*}}\mathcal{O}_X(K_X-Q))=\hat{B_i}$, hence for a general $y\in X$ with 	
                                            		\begin{equation*}y\notin Z:=\bigcap\limits_{L \in V^1(\mathop{f_{B_i*}}\mathcal{O}_X(K_X-Q))^c} Bs(\abs{K_X-Q+L}),	\end{equation*}
                                            		the set
                                            		
                                            		\begin{equation*}
                                            		\mathcal{B}_{f_{B_i}}(y)=\{L\in\hat{B_i}\vert \ y\notin Bs(\abs{K_X-Q+L})\}
                                            		\end{equation*}
                                            		
                                            		contains an open dense subset of $\hat{B_i}$ (cf. \cite[4.12]{blnp}), say $W$.
                                            		Taking the following groups into an evaluation map:
                                            		\begin{equation*} \bigoplus \limits_{R\in W}H^0(X,\mathcal{O}_X(2K_X+P)\otimes R^{-1})\otimes H^0(X,\mathcal{O}_X(K_X-Q+R))\rightarrow
                                            		H^0(3K_X+P-Q), \end{equation*} since $V+W=\hat{B_i}$, by (\ref{firsts}) it's easy to see that $H^0(3K_X+P-Q)$ has a section that passes through $x$ but does not pass through $y$. Therefore if we fix a general $P\in\mathrm{Pic}^0(X)_{tors}$, $\abs{3K_X+P-Q}$ separates $x$ and $y$, for any $x,y\in U_2\cap Z^c$ with $f_{B_i}(x)\ne f_{B_i}(y)$.
                                            		\medskip
                                            		
                                            		The last assertion now follows since for any $x,y\in X$, the equality $f_{B_i}(x)=f_{B_i}(y)$ for every $B_i$ implies that $d(x)=d(y)$, by Lemma \ref{tours}.
                                            	\end{proof}
                                            	\medskip
                                            	
                                            	Next we investigate the behavior under the fiber restriction map.
                                            	\medskip
                                            	\begin{lemma}\label{key}
                                            		Notations as in Theorem \ref{dt}, assume $B_i$ is a quotient abelian variety, $\F:=\F^j_{B_i}$ is a M-regular sheaf of rank 1 on its support, with the associated torsion line bundle $Q$. Take any $P\in \mathrm{Pic}^0(X)_{tors}$. Let $F$ be a general fiber of the composition $f_{B_i}:=p_{B_i}\circ a:X \rightarrow B:=f_{B_i}(X)$, then
                                            		the image of the natural restriction map: \begin{equation*}
                                            		H^0(X,3K_X+ P)\rightarrow H^0(F,3K_F+P)
                                            		\end{equation*} contains \begin{equation*}
                                            		\F\otimes \mathbb{C}_{f_{B_i}(F)}\subseteq \mathop{f_{B_i*}}(\mathcal{O}_X(K_X-Q))\otimes \mathbb{C}_{f_{B_i}(F)}\cong H^0(F,\mathcal{O}_F(K_F-Q))
                                            		\end{equation*} (via the natural base change isomorphism) considered as a subvector space of $H^0(F,\mathcal{O}_F(3K_F+P))$ by the inclusion induced by any nonzero section in the subvector space 
                                            		$ H^0(F,\mathcal{O}_F(2K_F+Q+P))$. (The inclusion $\F\hookrightarrow \mathop{f_{B_i*}}(\mathcal{O}_X(K_X-Q))$ is due to a version of the projection formula; we always identify $P, Q$ with its restriction to $F$.)
                                            	\end{lemma}
                                            	\begin{proof}
                                            		
                                            		We just prove the statement for $P=\mathcal{O}_X$, the same argument works for any $P\in \mathrm{Pic}^0(X)_{tors}$.

                                            		By construction $\mathop{f_{B_i*}}\mathcal{O}_X(2K_X+Q)$ is a M-regular sheaf on $B_i$, therefore continuously globally generated (Theorem \ref{Mtocgg}).
                                            		
                                            		This implies
                                            		\begin{equation}\label{firstcgg} \begin{split}
                                            		\bigoplus\limits_{R\in U}H^0(B_i,\mathop{f_{B_i*}}\mathcal{O}_X(2K_X+Q)\otimes R)\otimes R^{-1}\rightarrow \mathop{f_{B_i*}}\mathcal{O}_X(2K_X+Q)\otimes\mathbb{C}_b \end{split}
                                            		\end{equation}
                                            		is surjective for any $b\in B$ and any dense open $U\subseteq \mathrm{Pic}^0(B_i)$.
                                            		
                                            		The same argument for the M-regular sheaf $\mathcal{F}$ gives that $\mathcal{F}$ is continuously globally generated at any $b\in B_i$, and we have a surjection:
                                            		\begin{equation*}
                                            		\bigoplus\limits_{R\in \mathrm{Pic}^0(B_i)}H^0(\F\otimes R)\otimes R^{-1}\rightarrow \F\otimes\mathbb{C}_b.
                                            		\end{equation*}
                                            		As $\F$ is of rank one on its support, take a general $b\in f_{B_i}(X)$ we may assume $\F\otimes\mathbb{C}_b\cong \mathbb{C}$, and then there is a nonempty open $U_b\subseteq \mathrm{Pic}^0(B_i)$ such that 
                                            		\begin{equation}\label{scgg}
                                            		H^0(\F\otimes R)\otimes R^{-1}\rightarrow \F\otimes\mathbb{C}_b	
                                            		\end{equation} is surjective for any $R\in U_b$. 
                                            		
                                            		Take $U=-U_b$ in (\ref{firstcgg}) and compose it with (\ref{scgg}), we get a surjection
                                            		\begin{equation*} \begin{split}
                                            		\bigoplus\limits_{R\in U_b}H^0(B_i,\mathop{f_{B_i*}}\mathcal{O}_X(2K_X+Q)\otimes R^{-1})\otimes H^0(B_i, \F\otimes R) \rightarrow\\ \mathop{f_{B_i*}}\mathcal{O}_X(2K_X+Q)\otimes\mathbb{C}_b\cdot\F\otimes \mathbb{C}_b
                                            		\subseteq \mathop{f_{B_i*}}\mathcal{O}_X(3K_X)\otimes\mathbb{C}_b.
                                            		\end{split} 
                                            		\end{equation*}
                                            		Meanwhile, the inclusion and evaluation
                                            		\begin{equation*} \begin{split}
                                            		\bigoplus\limits_{R\in U_b}H^0(B_i,\mathop{f_{B_i*}}\mathcal{O}_X(2K_X+Q)\otimes R^{-1})\otimes H^0(B_i, \F\otimes R) \hookrightarrow \\
                                            		\bigoplus\limits_{R\in U_b}H^0(X,\mathcal{O}_X(2K_X+Q)\otimes R^{-1})\otimes H^0(X,\mathcal{O}_X(K_X-Q)\otimes R)\rightarrow	H^0(X,3K_X)\end{split} 
                                            		\end{equation*}
                                            		give the relation 
                                            		\begin{equation*}
                                            		\mathrm{Im}\{H^0(X,3K_X)\rightarrow H^0(X_F,3K_{X_F})\} \supseteq \mathop{f_{B_i*}}\mathcal{O}_X(2K_X+Q)\otimes\mathbb{C}_b\cdot\F\otimes \mathbb{C}_b	
                                            		\end{equation*}
                                            		for any general $b\in f_{B_i}(X)$ and $F:=f_{B_i}^{-1}(b)$. After base change this is nothing but the statement of the lemma.
                                            		
                                            	\end{proof}

                                            	We are now able to prove:
                                            	
                                            	\begin{theorem}\label{thmtors}Notations as above, if $X$ is of Albanese fiber genus 2
                                            		, then $\abs{3K_X+P}$ is birational, for general $P\in \mathrm{Pic}^0(X)_{tors}$. 
                                            	\end{theorem}
                                            	\begin{proof}
                                            		The proof is divided into five steps. In step 1, we show that $\abs{3K_X+P}$ is of degree at most 2 on a general  Albanese fiber. In step 2, we show that the Albanese of the morphism induced by $\abs{3K_X+P}$ is an isogeny. These two obstructions of birationality are reduced to two Galois morphisms in step 3. In step 4 and 5, we apply an argument of Jiang and Sun to show that both two Galois morphisms are in fact birational, which finishes the proof.
                                            		\medskip
                                            		
                                            		\textbf{Step 1.}		Take a general $P\in \mathrm{Pic}^0(X)_{tors}$ such that Lemma \ref{df} holds. 	Because $X$ is of general type, a general fiber of $a$ is a hyperelliptic curve of genus 2. Moreover, Theorem \ref{dt} tells us that for any $x\in X$, we have $a(x)+T\subseteq a(X)$.
                                            		
                                            		By the work of Jiang and Sun \cite{js}, we only need to consider the case that the cohomological loci is not the whole Picard zero group i.e. $V^0(K_X)\neq \mathrm{Pic}^0(A)$. In this case the sheaves $ \mathop{f_{B_i*}}\omega_X$ are torsion-free on the image of $X$, and by Theorem \ref{dt} the decomposition of $ \mathop{a_*}\omega_X$ consists of exactly two pull-back term $\F_i$ from the quotients $B_i$, with corresponding torsion elements $Q_i\in\mathrm{Pic}^0(X)_{tors}$, such that $\mathrm{rank}(\F_i\mid_{f_{B_i}(X)})=1$ ($i=1,2$). By Lemma \ref{tours} we have $\dim{B_i}\ge\dim f_{B_i}(X)>0$, since $f_{B_i}(X)$ generates $B_i$. (We remark that clearly $B_1\neq B_2$.)
                                            		
                                            		Because $\hat{B_i}$'s generate $\mathrm{Pic}^0(X)$ (Lemma \ref{tours}), there exist two elements $P_i\in \hat{B_i}\subseteq \mathrm{Pic}^0(X),\ i=1,2$ such that \begin{equation*}
                                            		P_1+Q_1=P_2+Q_2 .
                                            		\end{equation*}
                                            		
                                            		Let	$T_{B_i}$ be a general fiber of $q_i:A\rightarrow B_i$ over $f_{B_i}(X)$,
                                            		note that the sheaf $\mathop{q_i^*}\F_i$ is trivial restricted on $T_{B_i}$, and hence we have the natural surjections \begin{equation}\label{res}
                                            		(\F_i\otimes\mathbb{C}_{q_i(T_{B_i})}\cong H^0(T_{B_i}, \mathop{q_i^*}\F_i\mid_ {T_{B_i}}))\otimes \Om_{T_{B_i}}\longrightarrow \mathop{q_i^*}\F_i\mid_ {T_{B_i}}.
                                            		\end{equation}
                                            		
                                            		Take a general point $a_0\in {T_{B_i}}\subseteq a(X)$, let ${F_a}$ be the fiber of $a:X\rightarrow A$ over $a_0\in A$, and set $F_{B_i}:=a^{-1}(T_{B_i})$. By (\ref{res}) above, the image of the natural restriction map (via base change) 
                                            		\begin{equation*}H^0(F_{B_i},\omega_{F_{B_i}}\otimes Q_i^{-1})\supseteq\F_i\otimes\mathbb{C}_{q_i(T_{B_i})}\rightarrow H^0(F_a, \omega_{F_a}\otimes Q_i^{-1})\cong \mathop{a_*}\omega_X\otimes\mathbb{C}_{a_0}
                                            		\end{equation*} contains $\mathop{q_i^*}\F_i\otimes\mathbb{C}_{a_0}$, regarding the isomorphism $\bigoplus\limits_{i=1,2}\mathop{q_i^*}\F_i\otimes\mathbb{C}_{a_0}\cong\mathop{a_*}\omega_X\otimes\mathbb{C}_{a_0}$ induced by Theorem \ref{dt}.

                                            		Accordingly, the results above together with Lemma \ref{key} imply that the image of the restriction map of vector spaces $H^0(3K_X+P)\rightarrow H^0(3K_{F_a})$ contains the subvector space $s\cdot H^0(K_{F_a})$, where $s$ is a nonzero section in $H^0(2K_{F_a})$ coming from a section in $H^0(2K_X+P_1+Q_1+P)=H^0(2K_X+P_2+Q_2+P)$. Hence we have showed that $\abs{3K_X+P}$ is birational or of degree 2 on $F_a$ because so is $\abs{K_{F_a}}$.
                                            		\medskip
                                            		
                                            		\textbf{Step 2.}		
                                            		After some birational modifications, $\abs{3K_X+P}$ induces a map $\varphi:X\rightarrow Y$ onto a smooth variety $Y$ and hits into a commutative diagram:
                                            		\begin{equation*}
                                            		\begin{xy}
                                            		(0,20)*+{X}="v1";(20,20)*+{Y}="v2";(0,0)*+{X_A}="v3";(20,0)*+{W}="v4";(0,-20)*+{A}="v5";(20,-20)*+{A/{E}}="v6";(40,-20)*+{A/T}="v7";
                                            		{\ar@{->}_{\varphi}"v1";"v2"};
                                            		{\ar@{->}_{f}"v1";"v3"};
                                            		{\ar@{->}_{t}"v2";"v4"};
                                            		{\ar@{->}_{l}"v3";"v5"};
                                            		{\ar@{->}"v4";"v6"};
                                            		{\ar@{->}"v3";"v4"};
                                            		{\ar@/^1pc/^{b}"v2";"v6"};
                                            		{\ar@/_1pc/_{a}"v1";"v5"};
                                            		{\ar@{->}_{Alb(\varphi)}"v5";"v6"};
                                            		{\ar@{->}_{h}"v6";"v7"};
                                            		\end{xy}
                                            		\end{equation*}
                                            		
                                            		Where $a$ and $b$ are Albanese maps, $f$ and $t$ are their Stein factorizations. Let $G\subseteq {X\times Y\times A/T}$ be the graph of $id_X\times \varphi\times (h\circ b\circ\varphi):X\rightarrow X\times Y\times A/T$, and $p_{23}(G)$ be its projection in $Y\times A/T$. Then Lemma \ref{df} implies $p_{23}(G)\subseteq Y\times A/T$ is birational to $Y$ via the first projection, hence
                                            		$Alb(\varphi)$ is an isogeny over $A/T$ (with kernel $E$) by the universal property of the Albanese map.

                                            		\medskip
                                            		
                                            		\textbf{Step 3.}	We decompose $\varphi$ into two Galois morphisms. 
                                            		
                                            		We first assume that the linear system $\abs{3K_X+P}$ is of degree 2 on $F_a$, thus it induces the quotient of the hyperelliptic involution. The birational case could be treated similarly.
                                            		
                                            		Let $H$ be the graph of $X$ in $X\times Y\times X_A$, then the projection $g_1:H\rightarrow p_{23}(H)\subseteq Y\times X_A$ is of degree $2$ with Galois group $\mu_2$, and the further projection $g_2:p_{23}(H)\rightarrow Y$ is of degree $\sharp E$ (cardinality of $E$) with Galois group $E$, such that the composition $g_2\circ g_1$ is identified with $\varphi$ up to an isomorphism. Note that $E$ does not lift to a subgroup of $Bir(X)$ in general.

                                            		Let us do some birational modifications. 
                                            		
                                            		Let $V_0$ be a regularization of $p_{23}(H)$, i.e $V_0$ is a projective variety with an action of $E$ and a $E$-equivariant
                                            		birational map $\xi: V_0\dashrightarrow p_{23}(H)$ (see \cite[Theorem 3]{sum}). Taking an $E$-equivariant resolution
                                            		of singularities (see \cite{aw}), one can assume that $V_0$ is smooth. By the universal property of Albanese map we see $V_0\rightarrow A$ is $E$-equivariant, therefore $E$ acts freely on $V_0$. Replace $p_{23}(H)$ by $V_0$ and $Y$ by $V_0/E$, we may assume $g_2$ is an $E$-\'etale covering.

                                            		Doing the same process for $X$, we may assume $X$ is smooth and $\mu_2$ acts on $X$, with an $E$-equivariant morphism $X/\tau\rightarrow V_0$. Where $\tau$ is the generator of $\mu_2$, which is a biregular involution on $X$. After a further modification (cf. ~\cite[proof of Theorem 5.3, Step3]{js}), we assume that each connected component of the fixed loci of $\tau$ is a smooth divisor of $X$, thus $X/\tau$ is smooth. Replace $V_0$ by $V:=X/\tau$ and $V_0/E$ by $V/E$, we may assume $g_1$ is an ramified cover of degree 2, which is even flat (\cite[Corollary 18.17]{esb}):
                                            		\begin{equation*}
                                            		\begin{xy}
                                            		(0,20)*+{X}="v1";(20,20)*+{V:=X/\tau}="v2";(55,20)*+{V/E}="v3";(20,0)*+{V_0}="v4";(55,0)*+{V_0/E}="v5";(70,0)*+{Y}="v6";
                                            		{\ar@{->}_{g_1}"v1";"v2"};
                                            		{\ar@{->}_{g_2}"v2";"v3"};
                                            		{\ar@{->}_{}"v2";"v4"};
                                            		{\ar@{->}_{}"v4";"v5"};
                                            		{\ar@{->}_{}"v3";"v5"};
                                            		{\ar@{->}_{}"v5";"v6"};
                                            		{\ar@/^2pc/_{\varphi}"v1";"v3"};
                                            		\end{xy}
                                            		\end{equation*}

                                            		\medskip
                                            		
                                            		\textbf{Step 4.}

                                            		At this step, we show that $\abs{3K_X+P}$ induces a degree 1 map on $F_a$ (i.e. $g_1$ is birational), using an argument of ~\cite[Proof of Theorem 5.3]{js}.
                                            		
                                            		Assume the contrary. Since $g_1:X\rightarrow V$ is a finite Galois flat morphism, $\mathop{g_{1*}}\mathcal{O}_X$ is decomposed into a direct sum of rank one flat (invertible) eigensheaves. Then there is a divisor $D$ on $V$ with (\cite[2.44]{kol})
                                            		\begin{equation*}
                                            		\mathop{g_{1*}}\mathcal{O}_X=\mathcal{O}_V\oplus\mathcal{O}_V(-D).
                                            		\end{equation*}
                                            		
                                            		The multiplication $\mathcal{O}_V(-D)\otimes \mathcal{O}_V(-D)\rightarrow \mathcal{O}_V$ of the $\mathcal{O}_V$ algebra $\mathop{g_{1*}}\mathcal{O}_X$ gives an effective divisor $B\in\abs{2D}$.
                                            		By Stein factorization $X\cong \mathrm{\textbf{Spec}}_V(\mathop{g_{1*}}\mathcal{O}_X)$, hence $g_1$ is a double covering branched along $B$, and we have :
                                            		\begin{equation*}
                                            		\omega_X=\mathop{g_{1}^*}(\omega_V\otimes\mathcal{O}_V(D)).
                                            		\end{equation*}

                                            		Moreover, a general fiber $T'$ of $f'$ (the Stein factorization of $a'$, see the diagram below) is isomorphic to $\mathbb{P}^1$. Applying the Hurwitz formula for $F_a\rightarrow T'$, we have $\mathop{\mathrm{deg}}_T'(D)=g_{F_a}+1=3$.
                                            		
                                            		By the projection formula we write the decomposition of $\mathop{g_{1*}}(\omega_X^3\otimes P)$ :
                                            		\begin{equation*} \begin{split}
                                            		\mathop{g_{1*}}(\mathcal{O}_X(3K_X+P))=\mathcal{O}_V(3K_V+3D+P) \oplus
                                            		\mathcal{O}_V(3K_V+2D+P). \end{split}
                                            		\end{equation*}
                                            		Because $\mathop{\mathrm{deg}}_T(\mathcal{O}_T(3K_V+2D+P))\ge0$, both $\mathop{a'_*}(\mathcal{O}_V(3K_V+3D))$ and $\mathop{a'_*}(\mathcal{O}_V(3K_V+2D+P))$ are non-zero sheaves on $X_V$.
                                            		
                                            		On the other hand, since $\mathop{ a_*}(\mathcal{O}_X(3K_X+P))$ is a M-regular sheaf (Lemma \ref{genericM}), so are the both sheaves $\mathop{a'_*}(\mathcal{O}_V(3K_V+3D))$ and $\mathop{a'_*}(\mathcal{O}_V(3K_V+2D))$, which means \begin{equation*}H^0(\mathop{a'_*}(\mathcal{O}_V(3K_V+3D))\otimes Q)\ne 0\end{equation*} and \begin{equation*}H^0(\mathop{a'_*}(\mathcal{O}_V(3K_V+2D))\otimes Q)\ne 0\end{equation*} for all $Q\in \hat{A}$, by Lemma \ref{Mtocgg}. 
                                            		
                                            		But the non-zero group $H^0(\mathcal{O}_V(3K_V+2D+P))$ is the involution-anti-invariant part of $H^0(\mathop{g_{1*}}(\mathcal{O}_X(3K_X+P)))$, under the action of $\mu_2$, which means $H^0(X,\mathcal{O}_X(3K_X+P))$ has a section that separates $x\in F_a$ and its involution. Thus $\abs{3K_X+P}$ must induce a birational map on $F_a$.
                                            		
                                            		\begin{equation*}
                                            		\begin{xy}
                                            		(0,20)*+{X}="v1";(20,20)*+{V}="v2";(0,0)*+{X_A}="v3";(20,0)*+{X_V}="v4";(0,-20)*+{A}="v5";(20,-20)*+{A}="v6";(40,20)*+{Y}="v7";(40,0)*+{W}="v8";(40,-20)*+{A/E}="v9";
                                            		{\ar@{->}_{g_1}"v1";"v2"};
                                            		{\ar@{->}_{f}"v1";"v3"};
                                            		{\ar@{->}_{f'}"v2";"v4"};
                                            		{\ar@{->}"v3";"v5"};
                                            		{\ar@{->}"v4";"v6"};
                                            		{\ar@{->}^{\cong}"v3";"v4"};
                                            		{\ar@/^1pc/^{a'}"v2";"v6"};
                                            		{\ar@/_1pc/_{a}"v1";"v5"};
                                            		{\ar@{->}_{=}"v5";"v6"};
                                            		{\ar@{->}_{g_2}"v2";"v7"};
                                            		{\ar@{->}"v6";"v9"};
                                            		{\ar@/^1pc/^{b}"v7";"v9"};
                                            		{\ar@{->}_{t}"v7";"v8"};	
                                            		{\ar@{->}"v8";"v9"};
                                            		\end{xy}
                                            		\end{equation*}
                                            		
                                            		\medskip
                                            		
                                            		\textbf{Step 5.}			
                                            		We now use the same argument to gurantee the birationality of $\abs{3K_X+P}$.
                                            		
                                            		Since $\mathop{\mathrm{deg}}(g_1)=1$, we may assume $\varphi$ is identified with $g_2$, by the construction above.
                                            		
                                            		Since $\varphi$ is an \'etale covering, we write 
                                            		\begin{equation*}
                                            		\mathop{\varphi_*}\mathcal{O}_X=\bigoplus\limits_{i}\mathcal{O}_Y\otimes R_i.
                                            		\end{equation*}
                                            		Here $R_i\in \hat{E}\subseteq\hat{A/{E}},\ 1\le i\le\sharp{E}$ (by the dual exact sequence of abelian varieties) are torsion elements respect to the decomposition of eigensheaves under Galois morphism: $\mathop{Alb(\varphi)_*}\mathcal{O}_{A}=\bigoplus \mathcal{O}_{A/T_0}(R_i)$.
                                            		
                                            		Fix a $P_0\in \hat{A/E}$ with $\mathop{\varphi^*} P_0=P$. By the projection formula we write $\mathop{\varphi_*}(\omega_X^3\otimes P)$ :
                                            		\begin{equation*} \begin{split}
                                            		\mathop{\varphi_*}(\mathcal{O}_X(3K_X+P))=\bigoplus\limits_{i}\mathcal{O}_Y(3K_Y+P_0)\otimes R_i. \end{split}
                                            		\end{equation*}
                                            		Because $\mathop{\mathrm{deg}}_{F_t}(\mathcal{O}_{F_t}(3K_Y+P_0+R_i))>0$, all $\mathop{b_*}(\mathcal{O}_Y(3K_Y+P_0)\otimes R_i)$ are non-zero sheaves on $W$.
                                            		
                                            		Since $\mathop{(Alb(\varphi)\circ a)_*}(\mathcal{O}_X(3K_X+P))$ is a M-regular sheaf, so are the sheaves $\mathop{b_*}(\mathcal{O}_Y(3K_Y+P_0)\otimes R_i)$, which means 
                                            		\begin{equation*}
                                            		H^0(\mathop{b_*}(\mathcal{O}_Y(3K_Y+P_0)\otimes R_i)\otimes Q)\ne 0
                                            		\end{equation*}
                                            		for all $Q\in \hat {A/{E}}$. 
                                            		
                                            		Because the non-zero groups $H^0(\mathcal{O}_Y(3K_Y+3D+P_0)\otimes\mathcal{O}_Y(R_i))$ are the eigenspaces respect to the action of $E$ in $H^0(X,3K_X+P)=H^0(Y,\mathop{\varphi_*}(\omega_X^3\otimes P))$, $H^0(X,\mathcal{O}_X(3K_X+P))$ has sections that separate a general orbit of $E$ i.e. $\mathrm{deg(g_2)=1}$ and $E=\{0\}$. Thus $\abs{3K_X+P}$ is birational.

                                            	\end{proof}
                                            	\medskip
                                            	\noindent\textit{Proof of Theorem \ref{main}:}
                                            	Let $p_1,p_2$ be the first and second projection of $X\times \mathrm{Pic}^0(X)$, denote by $\F:=(\mathop{p_1^*}\omega_X^3)\otimes \Pd$, where $\Pd$ is the Poincar\'e bundle on $X\times \mathrm {Pic}^0(X)$. Then the surjection of the counit map $\mu: \mathop{p_2^*}\mathop{p_{2*}}\F\rightarrow \F$ onto its image induces a rational map (cf.\cite[ChII, Prop7.12]{h}):
                                            	\begin{equation*}
                                            	\Phi: X\times \mathrm{Pic}^0(X) \dashrightarrow \textbf{Proj}_{\mathrm{Pic}^0(X)}\bigoplus_{d>0} \mathrm{S}^d(\mathop{p_{2*}}\F).
                                            	\end{equation*}
                                            	
                                            	By generic flatness, there is a nonempty open $U_0\subseteq \mathrm{Pic}^0(X)$ over which $\F$ is flat, by Grauert's theorem there is a nonempty open $U\subseteq U_0$ over which $\mathop{p_{2*}}\F$ is locally free and has isomorphic base change map. We assume further that the map $\mu$ is surjective at the generic point of any fiber of $p_2$ over $U$. Then for any $P\in U$, the map $\Phi\vert_{X\times\{P\}}$ coincides with $\abs{3K_X+P}$ on a nonempty open set of $X$. Moreover, by Theorem \ref{js1} and the M-regularity of $\mathop{a_*}\omega_X^2\otimes Q$, we have $\dim\abs{3K_X+P}\ge 1$ for any $P\in \mathrm{Pic}^0(X)$.
                                            	
                                            	Let $W_0\subseteq X\times \mathrm{Pic}^0(X)$ be the open set of definition of $\Phi$. Then $W:=\Phi(W_0)$ is a constructible subset of $\textbf{Proj}_{\mathrm{Pic}^0(X)}\bigoplus_{d>0} \mathrm{S}^d(\mathop{p_{2*}}\F)$, which carries an integral scheme structure. Denote $\pi: W \rightarrow \mathrm{Pic}^0(X)$ the natural projection. Take a general fiber $W_P:=\pi^{-1}P$ over $P\in U$. As $\mathrm{Pic}^0(X)_{tors}$ is dense in $\mathrm{Pic}^0(X)$, by the semicontinuity of fiber dimension(\cite[ChII, Exercise 3.22]{h}) and Theorem \ref{thmtors} we have
                                            	\begin{equation*}
                                            	\dim W_P= \dim {W}_{P_0}= \dim X,
                                            	\end{equation*}
                                            	where $P_0\in\mathrm{Pic}^0(X)_{tors}$ is general. This implies $\dim W=2\dim X-1$ and $\Phi$ is generically finite.
                                            	
                                            	\medskip

                                            	By semicontinuity, generic base change (\cite[Proposition 4.1]{lps}) and Theorem \ref{thmtors} again, it's clear that the sheaf $\mathop{\Phi_*} \Om_{W_0}$ is of rank 1 on $W$, the statement follows.
                                            	\qed
                                            	
                                            	\medskip
                                            	\begin{remark}
                                            		When $X$ is of Albanese fiber genus $>2$, one can still use the argument above to show that, for a general $P\in\textrm{Pic}^0(X)$, $\abs{3K_X+P}$ is not generically finite of degree 2 (because degree 2 morphisms are always Galois).
                                            	\end{remark}
                                            	
                                            	At last, we give an example following the ideas of \cite[Example 1.13]{el} and \cite[Example 3.4]{js}.
                                            	
                                            	\begin{example}\label{examp}
                                            		(Ein-Lazarsfeld type threefold) Let $E$ be a elliptic curve and $C\rightarrow E$ a double covering by a curve $C$ of genus $\ge2$, and $\iota$ the corresponding involution on $C$. Let $\mathbb{P}^1\rightarrow \mathbb{P}^1$ be the 2nd power morphism with involution $\sigma_0$, $p:C'\rightarrow \mathbb{P}^1$ be a $\sigma_0$-equivariant double covering with the involution $\tau$ from a genus 2 curve $C'$, such that $\sigma_0$ lifts to another involution $\sigma$ on $C'$.
                                            		Take $X_0:=\frac{C\times C\times C'}{\iota\times\iota\times\tau}$, then there is an involution $\xi$ on $X_0$ induced by $\textrm{id}\times\iota\times\sigma$. 
                                            		Consider $X_1:=X_0/\xi$ into a diagram
                                            		\begin{equation*}
                                            		C\times C\times C'\xrightarrow{f} X_0 \xrightarrow{g} X_1 \xrightarrow{} A:=E\times E.
                                            		\end{equation*} 
                                            		We note that both $f$ and $g$ are finite morphisms of degree 2 and unramified in codimension $1$, hence $X_1$ has only canonical singularities (cf. \cite[Theorem 3.21]{kol}). Let $X$ be a desingularization of $X_1$, then the natural morphism $a:X\rightarrow A$ is the Albanese map, hence $X$ is of Albanese fiber genus 2 and of general type. Moreover, 
                                            		\begin{equation*}
                                            		V^0(K_X)=(\hat{E}\times\{\hat{0}\})\bigcup (\{\hat{0}\}\times \hat{E}) .\end{equation*}
                                            	\end{example}

                                            	\section*{Acknowledgements}
                                            	
                                            	The author thanks greatly to Zhi Jiang who, besides giving this problem, pointed out several significant mistakes in the original draft of this paper and made a variety of useful suggestions. The author also thanks his advisor Meng Chen for constant support and encouragement.

                                            	\section{References}
                                            	
                                            	\bibliographystyle{plain}
                                            	\bibliography{ref.bib}
                                            	
                                            \end{document}